\documentclass[a4paper, 12pt]{article}
\usepackage{amssymb}

\usepackage{amsmath}

\begin{document}

\begin{center}
{\Huge Recurrence and ergodicity in unital $\ast$-algebras }

{\Huge \vspace{10 mm} }{\normalsize Rocco Duvenhage and Anton Str\"{o}h }

\emph{Department of Mathematics and Applied Mathematics}

\emph{University of Pretoria, 0002 Pretoria, South Africa}{\normalsize \ }
\end{center}

\noindent\emph{Abstract.} Results concerning recurrence and ergodicity are
proved in an abstract Hilbert space setting based on the proof of
Khintchine's recurrence theorem for sets, and on the Hilbert space
characterization of ergodicity. These results are carried over to a
non-commutative $\ast $-algebraic setting using the GNS-construction. This
generalizes the corresponding measure theoretic results, in particular a
variation of Khintchine's Theorem for ergodic systems, where the image of
one set overlaps with another set, instead of with itself.

\section{Introduction}

The inspiration for this paper is the following theorem of Khintchine dating
from 1934 (see [\textbf{4}] for a proof):

\bigskip\noindent\textsc{Khintchine}'\textsc{s Theorem. }\emph{Let }$%
(X,\Sigma,\mu)$\emph{\ be a probability space (that is to say, }$\mu $\emph{%
\ is a measure on a }$\sigma$\emph{-algebra }$\Sigma$\emph{\ of subsets of a
set }$X$\emph{, with }$\mu(X)=1$\emph{), and consider a mapping }$%
T:X\rightarrow X$\emph{\ such that }$T^{-1}(S)\in\Sigma$\emph{\ and }$%
\mu(T^{-1}(S))\leq\mu(S)$\emph{\ for all }$S\in\Sigma$\emph{. Then for any }$%
A\in\Sigma$\emph{\ and }$\varepsilon>0$\emph{, the set} 
\begin{equation*}
E=\left\{ k\in\mathbb{N}:\mu\left( A\cap T^{-k}(A)\right) >\mu
(A)^{2}-\varepsilon\right\}
\end{equation*}
\emph{is relatively dense in }$\mathbb{N}=\{1,2,3,...\}$\emph{.}

\bigskip We will call $(X,\Sigma,\mu,T)$, as given above, a \emph{measure
theoretic dynamical system}. Recall that the relatively denseness of $E$ in $%
\mathbb{N}$ means that there exists an $n\in\mathbb{N}$ such that $%
E\cap\{j,j+1,...,j+n-1\}$ is non-empty for every $j\in\mathbb{N}$.
Khintchine's Theorem is an example of a recurrence result. It tells us that
for every $k\in E$, the set $A$ contains a set $A\cap T^{-k}(A)$ of measure
larger than $\mu(A)^{2}-\varepsilon$ which is mapped back into $A$ by $T^{k}$%
.

A question that arises from Khintchine's Theorem is whether, given $%
A,B\in\Sigma$ and $\varepsilon>0$, the set 
\begin{equation*}
F=\left\{ k\in\mathbb{N}:\mu\left( A\cap T^{-k}(B)\right) >\mu
(A)\mu(B)-\varepsilon\right\}
\end{equation*}
is relatively dense in $\mathbb{N}$. This is clearly not true in general,
for example if $T$ is the identity and $A$, $B$ and $\varepsilon$ are chosen
such that $\mu(A)\mu(B)>\varepsilon$ while $A\cap B$ is empty, then $F$ is
empty. $T$ has to ``mix'' the measure space sufficiently for $F$ to be
non-empty. In [\textbf{5}] it is shown for the case where $%
\mu(T^{-1}(S))=\mu(S)$ for all $S\in\Sigma$, that if for every pair $%
A,B\in\Sigma$ of positive measure there exists some $k\in\mathbb{N}$ such
that $\mu\left( A\cap T^{-k}(B)\right) >0$, then the dynamical system is
ergodic. Ergodicity therefore seems like the natural concept to use when
considering the question posed above. This is indeed what we will do.

The notion of ergodicity originally developed as a way to characterize
systems in classical statistical mechanics for which the time mean and the
phase space mean of any observable are equal. For our purposes it will be
most convenient to define ergodicity of a measure theoretic dynamical system 
$(X,\Sigma ,\mu,T)$ as follows (refer to [\textbf{4}], for example): $%
(X,\Sigma,\mu,T)$ is called \emph{ergodic }if the fixed points of the linear
Hilbert space operator $U:L^{2}(\mu)\rightarrow L^{2}(\mu):f\mapsto f\circ T$
form a one-dimensional subspace of $L^{2}(\mu)$. (It is easy to verify that $%
U$ is well-defined on $L^{2}(\mu)$.)

As we shall see, the ideas we have discussed so far are not really measure
theoretic in nature. This is in large part due to the fact that the proof of
Khintchine's Theorem is essentially a Hilbert space proof using the Mean
Ergodic Theorem. This proof can for the most part be written purely in
Hilbert space terms, hence giving an abstract Hilbert space result. Along
with the Hilbert space characterization of ergodicity given above, this
means that a fair amount of ergodic theory can be done purely in an abstract
Hilbert space setting. This is the approach taken in Section 3, using the
Mean Ergodic Theorem as the basic tool.

Having built up some ergodic theory in abstract Hilbert spaces, nothing is
to stop us from applying the results to mathematical structures other than
measure theoretic dynamical systems. The mathematical structure we will
consider is much more general than measure theoretic dynamical systems and
can easily be motivated as follows: From a measure theoretic dynamical
system $(X,\Sigma,\mu,T)$ we obtain the unital $\ast$-algebra $%
B_{\infty}(\Sigma)$ of all bounded complex-valued measurable functions
defined on $X$, and two linear mappings 
\begin{equation*}
\varphi:B_{\infty}(\Sigma)\rightarrow\mathbb{C}:f\mapsto\int fd\mu
\end{equation*}
and 
\begin{equation}
\tau:B_{\infty}(\Sigma)\rightarrow B_{\infty}(\Sigma):f\mapsto f\circ T 
\tag{1}
\end{equation}
with the following properties: $\varphi(1)=1$, $\varphi(f^{\ast}f)\geq0$, $%
\tau(1)=1$ and $\varphi(\tau(f)^{\ast}\tau(f))\leq\varphi(f^{\ast}f)$ for
all $f\in B_{\infty}(\Sigma)$, where $f^{\ast}=\overline{f}$ defines the
involution on $B_{\infty}(\Sigma)$, making it a $\ast$-algebra. We can view
this abstractly by replacing $B_{\infty}(\Sigma)$ with any unital $\ast $%
-algebra and considering linear mappings $\varphi$ and $\tau$ on it with the
properties mentioned above. (A \textit{unital }$\ast$\textit{-algebra} $%
\frak{A}$ is an algebra with an involution, and a unit element denoted by $1$%
, that is to say $1A=A=A1$ for all $A\in\frak{A}$. We will only work with
the case of complex scalars.) The most obvious generalization this brings is
that the unital $\ast$-algebra need not be commutative, for example the
bounded linear operators on a Hilbert space. Also note that $\tau$ in (1) is
a $\ast$-homomorphism of $B_{\infty}(\Sigma)$, but we will not need this
property of $\tau$ in the abstract $\ast$-algebraic setting. We describe the 
$\ast$-algebraic setting in more detail in Section 2, and in Section 4 the
Hilbert space results are applied to this setting using the GNS-construction.

In Section 5 we obtain the measure theoretic results as a special case, and
also briefly discuss another special case, namely von Neumann algebras.

\section{$*$-dynamical systems and ergodicity}

By a \emph{state} on a unital $*$-algebra $\frak{A}$ we mean a linear
functional $\varphi$ on $\frak{A}$ which is positive (i.e. $\varphi
(A^{*}A)\geq0$ for all $A\in\frak{A}$) with $\varphi(1)=1$. Motivated by our
remarks in Section 1, we give the following definition:

\bigskip\noindent\textsc{Definition 2.1. }Let $\varphi$ be a state on a
unital $*$-algebra $\frak{A}$. Consider any linear function $\tau:\frak{A}%
\rightarrow\frak{A}$ such that 
\begin{equation*}
\tau(1)=1
\end{equation*}
and 
\begin{equation*}
\varphi\left( \tau(A)^{*}\tau(A)\right) \leq\varphi(A^{*}A)
\end{equation*}
for all $A\in\frak{A}$. Then we call $(\frak{A},\varphi,\tau)$ a $*$\emph{%
-dynamical system}.

\bigskip Let $L(V)$ denote the algebra of all linear operators $V\rightarrow
V$ on the vector space $V$.

\bigskip\noindent\textsc{Definition 2.2. }Let $\varphi$ be a state on a
unital $\ast$-algebra $\frak{A}$. A \emph{cyclic representation} of $(\frak{A%
},\varphi)$ is a triple $(\frak{G},\pi,\Omega)$, where $\frak{G}$ is an
inner product space, $\pi:\frak{A}\rightarrow L(\frak{G})$ is linear with $%
\pi (1)=1$, $\pi(AB)=\pi(A)\pi(B)$, $\Omega\in\frak{G}$, $\pi(\frak{A}%
)\Omega=\frak{G}$, and $\left\langle \pi(A)\Omega,\pi(B)\Omega\right\rangle
=\varphi(A^{\ast}B)$, for all $A,B\in\frak{A}$.

\bigskip A cyclic representation as in Definition 2.2 exists by the
GNS-construction (refer to [\textbf{1}] for example, where the construction
is performed for C*-algebras, but it also works for unital $\ast $%
-algebras). We will not need the property $\pi (AB)=\pi (A)\pi (B)$ in this
paper however. The term ``cyclic'' refers to the fact that $\pi (\frak{A}%
)\Omega =\frak{G}$. Note that 
\begin{equation}
\iota :\frak{A}\rightarrow \frak{G}:A\mapsto \pi (A)\Omega  \tag{2}
\end{equation}
is a linear surjection such that $\iota (1)=\Omega $, and that 
\begin{equation}
U_{0}:\frak{G}\rightarrow \frak{G}:\iota (A)\mapsto \iota (\tau (A))  \tag{3}
\end{equation}
is a well-defined linear operator with $\left\| U_{0}\right\| \leq 1$ for $%
\tau $ as in Definition 2.1, since $\left\| \iota (\tau (A))\right\|
^{2}=\varphi (\tau (A)^{\ast }\tau (A))\leq \varphi (A^{\ast }A)=\left\|
\iota (A)\right\| ^{2}$. We define a seminorm $\left\| \cdot \right\|
_{\varphi }$ on $\frak{A}$ by 
\begin{equation*}
\left\| A\right\| _{\varphi }=\sqrt{\varphi \left( A^{\ast }A\right) }%
=\left\| \iota (A)\right\|
\end{equation*}
for all $A\in \frak{A}$. We now want to define the concept of ergodicity for
a $\ast $-dynamical system.

\bigskip \noindent \textsc{Definition 2.3. }A $\ast $-dynamical system $(%
\frak{A},\varphi ,\tau )$ is called \emph{ergodic} if it has the following
property: For any sequence $(A_{n})$ in $\frak{A}$ such that $\left\| \tau
(A_{n})-A_{n}\right\| _{\varphi }\rightarrow 0$ and such that for any $%
\varepsilon >0$ there exists an $N\in \mathbb{N}$ for which $\left\|
A_{m}-A_{n}\right\| _{\varphi }\leq \varepsilon $ if $m>N$ and $n>N$, it
follows that $\left\| A_{n}-\alpha \right\| _{\varphi }\rightarrow 0$ for
some $\alpha \in \mathbb{C}$.

\bigskip In Section 4 we will give a simple example of an ergodic $\ast $%
-dynamical system whose $\ast$-algebra is non-commutative. Recall that for
any vectors $x$ and $y$ in a Hilbert space $\frak{H}$, we denote by $%
x\otimes y$ the bounded linear operator $\frak{H\rightarrow H}$ defined by $%
(x\otimes y)z=x\left\langle y,z\right\rangle $. The motivation for
Definition 2.3 is the following proposition:

\bigskip\noindent\textsc{Proposition 2.4. }\emph{Consider a }$\ast $\emph{%
-dynamical system }$(\frak{A},\varphi,\tau)$\emph{\ and let }$U_{0}$\emph{\
be given by (3) in terms of any cyclic representation of }$(\frak{A}%
,\varphi) $\emph{. Let }$U:\frak{H}\rightarrow\frak{H}$\emph{\ be the
bounded linear extension of }$U_{0}$\emph{\ to the completion }$\frak{H}$%
\emph{\ of }$\frak{G}$\emph{, and let }$P$\emph{\ be the projection of }$%
\frak{H}$\emph{\ onto the subspace of fixed points of }$U$\emph{. Then }$(%
\frak{A},\varphi,\tau)$ \emph{is ergodic if and only if }$P=\Omega
\otimes\Omega$\emph{, that is to say, if and only if the fixed points of }$U$%
\emph{\ form a one-dimensional subspace of }$\frak{H}$\emph{.\bigskip}

\noindent\emph{Proof. }Since $\left\| \Omega\right\| ^{2}=\varphi(1^{\ast
}1)=1$, we know that $\Omega\otimes\Omega$ is the projection of $\frak{H}$
onto the one-dimensional subspace $\mathbb{C}\Omega$. Also note that $%
U\Omega=\Omega$, since $\Omega=\iota(1)$, hence $\mathbb{C}\Omega\subset P%
\frak{H}$.

Suppose $(\frak{A},\varphi ,\tau )$ is ergodic and let $x$ be a fixed point
of $U$. Consider any sequence $(x_{n})$ in $\frak{G}$ such that $%
x_{n}\rightarrow x$, say $x_{n}=\iota (A_{n})$. Then $\left\| \tau
(A_{n})-A_{n}\right\| _{\varphi }=\left\| Ux_{n}-x_{n}\right\| \rightarrow 0$%
, since $U$ is continuous, while for any $\varepsilon >0$ there exists some $%
N$ for which $\left\| A_{m}-A_{n}\right\| _{\varphi }=\left\|
x_{m}-x_{n}\right\| <\varepsilon $ if $m>N$ and $n>N$. Since $(\frak{A}%
,\varphi ,\tau )$ is ergodic, it follows that $\left\| x_{n}-\iota (\alpha
)\right\| =\left\| A_{n}-\alpha \right\| _{\varphi }\rightarrow 0$ for some $%
\alpha \in \mathbb{C}$, but then $x=\iota (\alpha )=\alpha \Omega $.
Therefore $P\frak{H}=\mathbb{C}\Omega $ which means that $P=\Omega \otimes
\Omega $.

Conversely, suppose $P=\Omega \otimes \Omega $ and consider any sequence $%
(A_{n})$ in $\frak{A}$ such that $\left\| \tau (A_{n})-A_{n}\right\|
_{\varphi }\rightarrow 0$ and such that for any $\varepsilon >0$ there
exists some $N$ for which $\left\| A_{m}-A_{n}\right\| _{\varphi
}<\varepsilon $ if $m>N$ and $n>N$. Then $x_{n}=\iota (A_{n})$ is a Cauchy
sequence and hence convergent in $\frak{H}$, since $\left\|
x_{m}-x_{n}\right\| =\left\| A_{m}-A_{n}\right\| _{\varphi }$. Say $%
x_{n}\rightarrow x$, then $Ux_{n}\rightarrow Ux$ since $U$ is continuous.
Since $\left\| Ux_{n}-x_{n}\right\| =\left\| \tau (A_{n})-A_{n}\right\|
_{\varphi }\rightarrow 0$, it follows that $Ux_{n}\rightarrow x$, hence $%
Ux=x $. This means that $x\in P\frak{H}$ which implies that $x=\alpha \Omega 
$ for some $\alpha \in \mathbb{C}$. Therefore $\left\| A_{n}-\alpha \right\|
_{\varphi }=\left\| x_{n}-\alpha \Omega \right\| \rightarrow 0$, and so we
conclude that $(\frak{A},\varphi ,\tau )$ is ergodic.$\square $

\bigskip Proposition 2.4 tells us that Definition 2.3 includes the measure
theoretic definition as a special case. This can be seen as follows: From a
measure theoretic dynamical system $(X,\Sigma,\mu,T)$ we obtain the $*$%
-dynamical system $(B_{\infty}(\Sigma),\varphi,\tau)$, where $\varphi
(f)=\int fd\mu$ and $\tau(f)=f\circ T$ for all $f\in B_{\infty}(\Sigma)$. A
cyclic representation of $(B_{\infty}(\Sigma),\varphi,\tau)$ is $(\frak{G}%
,\pi,\Omega)$ with $\frak{G}=\left\{ [g]:g\in B_{\infty}(\Sigma)\right\} $, $%
\pi(f)[g]=[fg]$ for all $f,g\in B_{\infty}(\Sigma)$, and $\Omega=[1]$, where 
$[g]$ denotes the equivalence class of all measurable complex-valued
functions on the measure space that are almost everywhere equal to $g$. The
completion of $\frak{G}$ is $L^{2}(\mu)$, and $U$ in Proposition 2.4 is now
given by 
\begin{equation*}
Uf=f\circ T
\end{equation*}
for all $f\in L^{2}(\mu)$, where here we have dropped the $[\cdot]$
notation, as is standard for $L^{2}$-spaces ($f$ and $f\circ T$ now denote
equivalence classes of functions). Proposition 2.4 tells us that $%
(B_{\infty}(\Sigma),\varphi,\tau)$ is ergodic if and only if the fixed
points of $U$ form a one dimensional subspace of $L^{2}(\mu)$, in other
words if and only if $(X,\Sigma,\mu,T)$ is ergodic, as was mentioned in
Section 1.

Finally we remark that we use Definition 2.3 as the definition of
ergodicity, since it is formulated purely in terms of the objects $\frak{A}$%
, $\varphi$ and $\tau$ appearing in the $\ast$-dynamical system $(\frak{A}%
,\varphi,\tau)$, unlike Proposition 2.4 which involves a cyclic
representation of these objects. However, as a characterization of
ergodicity, Proposition 2.4 is generally easier to use. Of course, one might
wonder if Definition 2.3 could not be simplified by using a single element
rather than a sequence. With $U$ as in Proposition 2.4, and $x=\iota(A)$ for
some $A\in\frak{A}$, we have $Ux=x$ if and only if $\left\| Ux-x\right\| =0$%
, which is equivalent to $\left\| \tau(A)-A\right\| _{\varphi}=0$. For
ergodicity we need this to imply that $x=\alpha\Omega$ for some $\alpha\in%
\mathbb{C}$, which is equivalent to $\left\| A-\alpha\right\|
_{\varphi}=\left\| x-\alpha \Omega\right\| =0$. However, we cannot define
ergodicity as ``$\left\| \tau(A)-A\right\| _{\varphi}=0$ implies that $%
\left\| A-\alpha\right\| _{\varphi}=0$ for some $\alpha\in\mathbb{C}$'',
since Proposition 2.4 would no longer hold: There would be examples of
ergodic $\ast$-dynamical systems for which the fixed points of $U$ do not
form a one-dimensional subspace of $\frak{H}$. (In the Appendix we give such
an example.) Our theory would then fall apart, since much of our later work
is based on the fact that for ergodic systems the fixed point space of $U$
is one-dimensional. For example, the characterization of ergodicity in terms
of the equality of means of the sort mentioned in Section 1 (but extended to 
$\ast$-dynamical systems), implies this one-dimensionality. Also, this
one-dimensionality is used in our proof of the variation of Khintchine's
Theorem mentioned in Section 1. (See Sections 3 and 4 for details.) The use
of a sequence rather than a single element is therefore necessary in
Definition 2.3.

\section{Some ergodic theory in Hilbert spaces}

Our main tool in this section is the

\bigskip\noindent\textsc{Mean Ergodic Theorem. }\emph{Consider a linear
operator }$U:\frak{H}\rightarrow\frak{H}$\emph{\ with }$\left\| U\right\|
\leq1$\emph{\ on a Hilbert space }$\frak{H}$\emph{\ . Let }$P$\emph{\ be the
projection of }$\frak{H}$\emph{\ onto the subspace of fixed points of }$U $%
\emph{. For any }$x\in\frak{H}$\emph{\ we then have } 
\begin{equation*}
\frac{1}{n}\sum_{k=0}^{n-1}U^{k}x\rightarrow Px
\end{equation*}
\emph{as }$n\rightarrow\infty$\emph{.}

\bigskip Refer to [\textbf{4}] for a proof. We now state and prove a
generalized Hilbert space version of Khintchine's Theorem:

\bigskip\noindent\noindent\textsc{Theorem 3.1. }\emph{Let }$\frak{H}$\emph{, 
}$U$\emph{\ and }$P$\emph{\ be as in the Mean Ergodic Theorem above.
Consider any }$x,y\in\frak{H}$\emph{\ and }$\varepsilon>0$\emph{. Then the
set } 
\begin{equation*}
E=\left\{ k\in\mathbb{N}:\left| \left\langle x,U^{k}y\right\rangle \right|
>\left| \left\langle x,Py\right\rangle \right| -\varepsilon\right\}
\end{equation*}
\emph{is relatively dense in }$\mathbb{N}$\emph{.}

\bigskip\noindent\emph{Proof. }The proof is essentially the same as that of
Khintchine's Theorem. By the Mean Ergodic Theorem there exists an $n\in%
\mathbb{N}$ such that 
\begin{equation*}
\left\| \frac{1}{n}\sum_{k=0}^{n-1}U^{k}y-Py\right\| <\frac{\varepsilon }{%
\left\| x\right\| +1}\text{.}
\end{equation*}
Since $UPy=Py$ and $\left\| U\right\| \leq1$, it follows for any $j\in%
\mathbb{N}$ that 
\begin{equation*}
\left\| \frac{1}{n}\sum_{k=j}^{j+n-1}U^{k}y-Py\right\| \leq\left\| \frac
{1}{n}\sum_{k=0}^{n-1}U^{k}y-Py\right\| <\frac{\varepsilon}{\left\|
x\right\| +1}
\end{equation*}
and therefore 
\begin{equation*}
\left| \left\langle x,\frac{1}{n}\sum_{k=j}^{j+n-1}U^{k}y-Py\right\rangle
\right| \leq\left\| x\right\| \left\| \frac{1}{n}\sum_{k=j}^{j+n-1}U^{k}y-Py%
\right\| <\varepsilon\text{.}
\end{equation*}
Hence 
\begin{equation*}
\left| \left\langle x,Py\right\rangle \right| -\varepsilon<\left| \frac
{1}{n}\sum_{k=j}^{j+n-1}\left\langle x,U^{k}y\right\rangle \right| \leq 
\frac{1}{n}\sum_{k=j}^{j+n-1}\left| \left\langle x,U^{k}y\right\rangle
\right|
\end{equation*}
and so $\left| \left\langle x,U^{k}y\right\rangle \right| >\left|
\left\langle x,Py\right\rangle \right| -\varepsilon$ for some $k\in
\{j,j+1,...,j+n-1\}$, in other words $E$ is relatively dense in $\mathbb{N}$.%
$\square$

\bigskip Khintchine's Theorem corresponds to the case where $y=x$. The
following two propositions are the Hilbert space building blocks for two
characterizations of ergodicity to be considered in the next section.

\bigskip\noindent\textsc{Proposition 3.2. }\emph{Let }$\frak{H}$\emph{, }$U $%
\emph{\ and }$P$\emph{\ be as in the Mean Ergodic Theorem above. Consider an 
}$\Omega\in\frak{H}$\emph{\ and let }$\frak{T}$\emph{\ be any total set in }$%
\frak{H}$\emph{. Then the following hold:}

\emph{(i) If }$P=\Omega\otimes\Omega$\emph{, then} 
\begin{equation}
\left\| \frac{1}{n}\sum_{k=0}^{n-1}U^{k}y-\Omega\left\langle \Omega
,y\right\rangle \right\| \rightarrow0  \tag{4}
\end{equation}
\emph{as }$n\rightarrow\infty$\emph{, for every }$y\in\frak{H}$\emph{.}

\emph{(ii) If (4) holds for every }$y\in\frak{T}$\emph{, then }$P=\Omega
\otimes\Omega$\emph{.}

\bigskip\noindent\emph{Proof. }By the Mean Ergodic Theorem we know that 
\begin{equation}
\left\| \frac{1}{n}\sum_{k=0}^{n-1}U^{k}y-Py\right\| \rightarrow0  \tag{5}
\end{equation}
for every $y\in\frak{H}$ as $n\rightarrow\infty$, but for $P=\Omega
\otimes\Omega$ we have $Py=\Omega\left\langle \Omega,y\right\rangle $ and
this proves (i).

To prove (ii), consider any $y\in\frak{T}$. From (4) and (5) it then follows
that $Py=\Omega\left\langle \Omega,y\right\rangle =(\Omega\otimes\Omega)y$.
Since by definition the linear span of $\frak{T}$ is dense in $\frak{H}$,
and since $P$ and $\Omega\otimes\Omega$ are bounded (and hence continuous)
linear operators on $\frak{H}$, we conclude that $P=\Omega\otimes\Omega$.$%
\square$

\bigskip\noindent\textsc{Proposition 3.3. }\emph{Let }$\frak{H}$\emph{, }$U $%
\emph{\ and }$P$\emph{\ be as in the Mean Ergodic Theorem above. Consider an 
}$\Omega\in\frak{H}$\emph{\ and let }$\frak{S}$\emph{\ and }$\frak{T}$\emph{%
\ be total sets in }$\frak{H}$\emph{. Then the following hold:}

\emph{(i) If }$P=\Omega\otimes\Omega$\emph{, then} 
\begin{equation}
\frac{1}{n}\sum_{k=0}^{n-1}\left\langle x,U^{k}y\right\rangle \rightarrow
\left\langle x,\Omega\right\rangle \left\langle \Omega,y\right\rangle 
\tag{6}
\end{equation}
\emph{as }$n\rightarrow\infty$\emph{, for all }$x,y\in\frak{H}$\emph{.}

\emph{(ii) If (6) holds for all }$x\in\frak{S}$\emph{\ and }$y\in\frak{T}$%
\emph{, then }$P=\Omega\otimes\Omega$\emph{.}

\bigskip\noindent\emph{Proof. }Statement (i) follows immediately from
Proposition 3.2(i) by simply taking the inner product of $x$ with the
expression inside the norm in (4).

To prove (ii), consider any $x\in\frak{S}$ and $y\in\frak{T}$. From the Mean
Ergodic Theorem it follows that 
\begin{equation*}
\frac{1}{n}\sum_{k=0}^{n-1}\left\langle x,U^{k}y\right\rangle \rightarrow
\left\langle x,Py\right\rangle
\end{equation*}
as $n\rightarrow\infty$. Combining this with (6) we see that $\left\langle
x,Py\right\rangle =\left\langle x,\Omega\right\rangle \left\langle
\Omega,y\right\rangle =\left\langle x,(\Omega\otimes\Omega)y\right\rangle $.
Since the linear span of $\frak{S}$ is dense in $\frak{H}$, this implies
that $Py=(\Omega\otimes\Omega)y$. Hence $P=\Omega\otimes\Omega$ as in the
proof of Proposition 3.2(ii).$\square$

\bigskip The reason for using total sets will become clear in Sections 4 and
5.

\section{Ergodic results for $*$-dynamical systems}

In this section we carry the results of Section 3 over to $*$-dynamical
systems using cyclic representations. Firstly we give a $*$-dynamical
generalization of Khintchine's Theorem which follows from Theorem 3.1:

\bigskip\noindent\textsc{Theorem 4.1. }\emph{Let }$(\frak{A},\varphi,\tau )$%
\emph{\ be a }$*$\emph{-dynamical system, and consider any }$A\in\frak{A}$%
\emph{\ and }$\varepsilon>0$\emph{. Then the set } 
\begin{equation*}
E=\left\{ k\in\mathbb{N}:\left| \varphi\left( A^{*}\tau^{k}(A)\right)
\right| >\left| \varphi(A)\right| ^{2}-\varepsilon\right\}
\end{equation*}
\emph{is relatively dense in }$\mathbb{N}$\emph{.}

\bigskip\noindent\emph{Proof. }Let $U$ and $P$ be defined as in Proposition
2.4 in terms of any cyclic representation of $(\frak{A},\varphi) $. Set $%
x=\iota(A)$. From (3) it is clear that $\Omega=\iota(1)$ is a fixed point of 
$U$, so $\left\langle \Omega,x\right\rangle =\left\langle P\Omega
,x\right\rangle =\left\langle \Omega,Px\right\rangle $. It follows that $%
\left| \varphi(A)\right| =\left| \varphi(1^{*}A)\right| =\left| \left\langle
\Omega,x\right\rangle \right| \leq\left\| \Omega\right\| \left\| Px\right\|
=\left\| Px\right\| $. We also have $\varphi(A^{*}\tau^{k}(A))=\left\langle
x,U^{k}x\right\rangle $. Hence by Theorem 3.1, with $y=x$, the set $E$ is
relatively dense in $\mathbb{N}$.$\square$

\bigskip A C*-algebraic version of Theorem 4.1 was previously obtained in [%
\textbf{3}]. Next we use Theorem 3.1 to prove a variant of Theorem 4.1:

\bigskip\noindent\textsc{Theorem 4.2. }\emph{Let }$(\frak{A},\varphi,\tau )$%
\emph{\ be an ergodic }$*$\emph{-dynamical system, and consider any }$A,B\in%
\frak{A}$\emph{\ and }$\varepsilon>0$\emph{. Then the set } 
\begin{equation*}
E=\left\{ k\in\mathbb{N}:\left| \varphi\left( A\tau^{k}(B)\right) \right|
>\left| \varphi(A)\varphi(B)\right| -\varepsilon\right\}
\end{equation*}
\emph{is relatively dense in }$\mathbb{N}$\emph{.}

\bigskip \noindent \emph{Proof. }Let $U$ and $P$ be defined as in
Proposition 2.4 in terms of any cyclic representation of $(\frak{A},\varphi
) $. Set $x=\iota (A^{\ast })$ and $y=\iota (B)$. By Proposition 2.4 we have 
$Px=\alpha \Omega $ and $Py=\beta \Omega $ where $\overline{\alpha }%
=\left\langle x,\Omega \right\rangle =\varphi (A^{\ast \ast }1)=\varphi (A)$
and $\beta =\varphi (B)$. Therefore $\left| \left\langle x,Py\right\rangle
\right| =\left| \left\langle Px,Py\right\rangle \right| =\left| \overline{%
\alpha }\beta \right| \left\| \Omega \right\| ^{2}=\left| \varphi (A)\varphi
(B)\right| $. Furthermore, $\varphi (A\tau ^{k}(B))=\left\langle
x,U^{k}y\right\rangle $. Hence $E$ is relatively dense in $\mathbb{N}$ by
Theorem 3.1.$\square $

\bigskip We are now going to prove two characterizations of ergodicity using
Propositions 3.2 and 3.3 respectively. But first we need to consider a
notion of totality of a set in a unital $\ast$-algebra. (Remember that an
abstract unital $\ast$-algebra has no norm.)

\bigskip\noindent\textsc{Definition 4.3. }Let $\varphi$ be a state on a
unital $*$-algebra $\frak{A}$. A subset $\frak{T}$ of $\frak{A}$ is called $%
\varphi $\emph{-dense} in $\frak{A}$ if it is dense in the seminormed space $%
(\frak{A},\left\| \cdot\right\| _{\varphi})$. A subset $\frak{T}$ of $\frak{A%
}$ is called $\varphi$\emph{-total} in $\frak{A}$ if the linear span of $%
\frak{T}$ is $\varphi$-dense in $\frak{A}$.

\bigskip Trivially, a unital $*$-algebra is $\varphi$-total in itself for
any state $\varphi$.

\bigskip\noindent\textsc{Lemma 4.4. }\emph{Let }$\varphi$\emph{\ be a state
on a unital }$*$\emph{-algebra }$\frak{A}$\emph{, and consider any subset }$%
\frak{T}$\emph{\ of }$\frak{A}$\emph{. Let }$\iota$ \emph{be given by (2) in
terms of any cyclic representation of }$(\frak{A},\varphi)$\emph{, and let }$%
\frak{H}$\emph{\ be the completion of }$\frak{G}$\emph{. Then }$\frak{T}$ 
\emph{is }$\varphi$\emph{-total in }$\frak{A}$\emph{\ if and only if} $\iota(%
\frak{T})$\emph{\ is total in }$\frak{H}$\emph{.}

\bigskip\noindent\emph{Proof. }Suppose $\frak{T}$ is $\varphi$-total in $%
\frak{A}$, that is to say the linear span $\frak{B}$ of $\frak{T}$ is $%
\varphi$-dense in $\frak{A}$. Then $\iota(\frak{B})$ is dense in $\frak{G}%
=\iota(\frak{A})$, since for any $A\in\frak{A}$ there exists a sequence $%
(A_{n})$ in $\frak{B}$ such that $\left\| \iota(A_{n})-\iota(A)\right\|
=\left\| A_{n}-A\right\| _{\varphi}\rightarrow0$. But by definition $\frak{G}
$ is dense in $\frak{H}$, hence $\iota(\frak{B})$ is dense in $\frak{H}$.
Since $\iota$ is linear, this means that $\iota(\frak{T})$ is total in $%
\frak{H}$.

Conversely, suppose $\iota(\frak{T})$ is total in $\frak{H}$, then $\iota(%
\frak{B})$ is dense in $\frak{H}$. It follows that $\frak{B}$ is $\varphi$%
-dense in $\frak{A}$, since for any $A\in\frak{A}$ there exists a sequence $%
(A_{n})$ in $\frak{B}$ such that $\left\| A_{n}-A\right\| _{\varphi}=\left\|
\iota(A_{n})-\iota(A)\right\| \rightarrow0$. In other words, $\frak{T}$ is $%
\varphi$-total in $\frak{A}$.$\square$

\bigskip\noindent\textsc{Proposition 4.5. }\emph{Let }$(\frak{A},\varphi
,\tau)$\emph{\ be a }$*$\emph{-dynamical system, and consider any }$\varphi $%
\emph{-total set }$\frak{T}$\emph{\ in }$\frak{A}$\emph{. Then the following
hold:}

\emph{(i) If }$(\frak{A},\varphi,\tau)$\emph{\ is ergodic, then} 
\begin{equation}
\left\| \frac{1}{n}\sum_{k=0}^{n-1}\tau^{k}(A)-\varphi(A)\right\| _{\varphi
}\rightarrow0  \tag{7}
\end{equation}
\emph{as }$n\rightarrow\infty$\emph{, for every }$A\in\frak{A}$\emph{.}

\emph{(ii) If (7) holds for every }$A\in\frak{T}$\emph{, then }$(\frak{A}%
,\varphi,\tau)$\emph{\ is ergodic.}

\bigskip \noindent \emph{Proof. }Let $U$ and $P$ be defined as in
Proposition 2.4 in terms of any cyclic representation of $(\frak{A},\varphi
) $. Suppose $(\frak{A},\varphi ,\tau )$ is ergodic. For any $A\in \frak{A}$
we then have 
\begin{equation}
\left\| \frac{1}{n}\sum_{k=0}^{n-1}\tau ^{k}(A)-\varphi (A)\right\|
_{\varphi }=\left\| \frac{1}{n}\sum_{k=0}^{n-1}U^{k}\iota (A)-\iota \left(
\varphi (A)\right) \right\| \rightarrow 0  \tag{8}
\end{equation}
as $n\rightarrow \infty $, by Proposition 3.2(i) and Proposition 2.4, since $%
\iota \left( \varphi (A)\right) =\iota (1)\varphi (A)=\Omega \varphi
(1^{\ast }A)=\Omega \left\langle \Omega ,\iota (A)\right\rangle $. This
proves (i).

Now suppose (7), and therefore (8), hold for every $A\in \frak{T}$. Since $%
\iota (\frak{T})$ is total in $\frak{H}$ according to Lemma 4.4, it follows
from Proposition 3.2(ii) and the identity $\iota \left( \varphi (A)\right)
=\Omega \left\langle \Omega ,\iota (A)\right\rangle $, that $P=\Omega
\otimes \Omega $. So $(\frak{A},\varphi ,\tau )$ is ergodic by Proposition
2.4, confirming (ii).$\square $

\bigskip In the spirit of the original motivation behind the concept of
ergodicity, this proposition characterizes ergodic $\ast$-dynamical systems
as those for which the \emph{time mean} of each element $A$ of the $\ast$%
-algebra converges in the seminorm $\left\| \cdot\right\| _{\varphi}$ to the
``phase space'' mean $\varphi(A)$. A better name for the latter would be the 
\emph{system mean} in this case, since there is no phase space involved. For
a measure theoretic dynamical system $(X,\Sigma,\tau,\mu)$, the state $%
\varphi$ is given by $\varphi(f)=\int fd\mu$ which is indeed the phase space
mean of $f\in B_{\infty}(\Sigma)$, where $X$ is the phase space. We will
come back to this in Section 5.

For any subset $\frak{S}$ of a $*$-algebra, we write $\frak{S}%
^{*}=\{A^{*}:A\in\frak{S}\}$.

\bigskip\noindent\textsc{Proposition 4.6. }\emph{Let }$(\frak{A},\varphi
,\tau)$\emph{\ be a }$*$\emph{-dynamical system, and consider any }$\varphi $%
\emph{-total sets $\frak{S}$ and }$\frak{T}$\emph{\ in }$\frak{A}$\emph{.
Then the following hold:}

\emph{(i) If }$(\frak{A},\varphi,\tau)$\emph{\ is ergodic, then} 
\begin{equation}
\frac{1}{n}\sum_{k=0}^{n-1}\varphi\left( A\tau^{k}(B)\right) \rightarrow
\varphi(A)\varphi(B)  \tag{9}
\end{equation}
\emph{as }$n\rightarrow\infty$\emph{, for all }$A,B\in\frak{A}$\emph{.}

\emph{(ii) If (9) holds for all }$A\in\frak{S}^{*}$\emph{\ and }$B\in\frak{T}
$\emph{, then }$(\frak{A},\varphi,\tau)$\emph{\ is ergodic.}

\bigskip\noindent\emph{Proof. }Let $U$ and $P$ be defined as in Proposition
2.4 in terms of any cyclic representation of $(\frak{A},\varphi) $. Suppose $%
(\frak{A},\varphi,\tau)$ is ergodic. Then $P=\Omega\otimes\Omega$ by
Proposition 2.4, and so by Proposition 3.3(i) it follows that 
\begin{equation}
\frac{1}{n}\sum_{k=0}^{n-1}\varphi\left( A\tau^{k}(B)\right) =\frac{1}{n}%
\sum_{k=0}^{n-1}\left\langle \iota(A^{*}),U^{k}\iota(B)\right\rangle
\rightarrow\varphi(A)\varphi(B)  \tag{10}
\end{equation}
as $n\rightarrow\infty$, since $\left\langle
\iota(A^{*}),\Omega\right\rangle =\varphi(A)$ and $\left\langle
\Omega,\iota(B)\right\rangle =\varphi(B)$, as in the proof of Theorem 4.2.
This proves (i). (Alternatively, (i) can be derived from Proposition 4.5(i)
using the Cauchy-Schwarz inequality $\left| \varphi(AC)\right| \leq\left\|
A^{*}\right\| _{\varphi}\left\| C\right\| _{\varphi}$ with $C=\frac{1}{n}%
\sum_{k=0}^{n-1}\tau^{k}(B)-\varphi(B)$. This is essentially how Proposition
3.3(i) was derived from Proposition 3.2(i).)

Now suppose (9), and therefore (10), hold for all $A\in \frak{S}^{\ast }$
and $B\in \frak{T}$. Since $\iota (\frak{S})$ and $\iota (\frak{T})$ are
total in $\frak{H}$ according to Lemma 4.4, it follows from Proposition
3.3(ii) and the identities $\left\langle \iota (A^{\ast }),\Omega
\right\rangle =\varphi (A)$ and $\left\langle \Omega ,\iota (B)\right\rangle
=\varphi (B)$, that $P=\Omega \otimes \Omega $. So $(\frak{A},\varphi ,\tau
) $ is ergodic by Proposition 2.4, confirming (ii).$\square $

\bigskip This characterizes ergodicity in terms of \emph{mixing}. We now
give a simple example of an ergodic $*$-dynamical system whose $*$-algebra
is non-commutative:

\bigskip\noindent\textsc{Example 4.7. }Let $\frak{A}$ be the unital $*$%
-algebra of $2\times2$-matrices with entries in $\mathbb{C}$, the involution
being the conjugate transpose. Let $\varphi$ be the normalized trace on $%
\frak{A}$, that is to say $\varphi=\frac{1}{2}$Tr. Define $\tau:\frak{A}%
\rightarrow\frak{A}$ by 
\begin{equation*}
\tau\left( 
\begin{array}{ll}
a_{11} & a_{12} \\ 
a_{21} & a_{22}
\end{array}
\right) =\left( 
\begin{array}{cc}
a_{22} & c_{1}a_{12} \\ 
c_{2}a_{21} & a_{11}
\end{array}
\right)
\end{equation*}
for some fixed $c_{1},c_{2}\in\mathbb{C}$ with $\left| c_{1}\right| \leq1$, $%
\left| c_{2}\right| \leq1$, $c_{1}\neq1$ and $c_{2}\neq1$. The conditions $%
\left| c_{1}\right| \leq1$ and $\left| c_{2}\right| \leq1$ are necessary and
sufficient for $(\frak{A},\varphi,\tau)$ to be a $*$-dynamical system. Note
that for any $c\in\mathbb{C}$ with $\left| c\right| \leq1$, it follows from
the Mean Ergodic Theorem that 
\begin{equation*}
\frac{1}{n}\sum_{k=0}^{n-1}c^{k}
\end{equation*}
converges to $0$ if $c\neq1$, and to $1$ otherwise. Using this fact and
Proposition 4.6(ii) with $\frak{S}=\frak{T}=\frak{A}$ (and some
calculations), it can be verified that the conditions $c_{1}\neq1$ and $%
c_{2}\neq1$ are necessary and sufficient for $(\frak{A},\varphi,\tau)$ to be
ergodic, assuming that $\left| c_{1}\right| \leq1$ and $\left| c_{2}\right|
\leq1$.

\section{Measure theory and von Neumann algebras}

As was mentioned in Section 2, from a measure theoretic dynamical system $%
(X,\Sigma,\mu,T)$ we obtain the $*$-dynamical system $(B_{\infty}(\Sigma),%
\varphi,\tau)$, where $\varphi(f)=\int fd\mu$ and $\tau(f)=f\circ T$. This
allows us to apply the results of Section 4 to measure theoretic dynamical
systems. For example, if $(X,\Sigma,\mu,T)$ is ergodic, then we know from
Section 2 that $(B_{\infty}(\Sigma),\varphi,\tau)$ is ergodic. Hence for
this $*$-dynamical system Theorem 4.2 tells us that for any $A,B\in\Sigma$
and $\varepsilon>0$, the set 
\begin{equation*}
\left\{ k\in\mathbb{N}:\left| \varphi\left( \chi_{A}\tau^{k}(\chi
_{B})\right) \right| >\left| \varphi(\chi_{A})\varphi(\chi_{B})\right|
-\varepsilon\right\}
\end{equation*}
is relatively dense in $\mathbb{N}$, but this set is exactly the set $F$
from Section 1. (Here $\chi$ denotes characteristic functions.) So we have
answered our original question:

\bigskip\noindent\textsc{Corollary 5.1. }\emph{Let }$(X,\Sigma,\mu ,T)$\emph{%
\ be an ergodic measure theoretic dynamical system. Then for any }$%
A,B\in\Sigma$\emph{\ and }$\varepsilon>0$\emph{, the set} 
\begin{equation*}
F=\left\{ k\in\mathbb{N}:\mu\left( A\cap T^{-k}(B)\right) >\mu
(A)\mu(B)-\varepsilon\right\}
\end{equation*}
\emph{is relatively dense in }$\mathbb{N}$\emph{.}

\bigskip This result says that for every $k\in F$, the set $A$ contains a
set $A\cap T^{-k}(B)$ of measure larger than $\mu(A)\mu(B)-\varepsilon$,
which is mapped into $B$ by $T^{k}$. Using a similar argument, Khintchine's
Theorem follows from Theorem 4.1.

Likewise, Propositions 4.5 and 4.6 can be applied to the measure theoretic
case. For example, Proposition 4.5(i) tells us that if $(X,\Sigma,\mu,T)$ is
ergodic, then 
\begin{equation}
\int\left| \frac{1}{n}\sum_{k=0}^{n-1}f\circ T^{k}-\varphi(f)\right|
^{2}d\mu\rightarrow0  \tag{11}
\end{equation}
as $n\rightarrow\infty$, for every $f\in B_{\infty}(\Sigma)$. Note that this
result is not pointwise and is therefore not quite as strong as the usual
measure theoretic statement of equality of the time mean and the phase space
mean. This is of course where Birkhoff's Pointwise Ergodic Theorem comes
into play (see for example [\textbf{4}]).

What about the converse? Well, in order to effectively apply Propositions
4.5(ii) and 4.6(ii) to the measure theoretic case, we need to know what the
measure theoretic significance of a $\varphi$-total set in $%
B_{\infty}(\Sigma)$ is. The basic fact we will use is the following simple
proposition which follows from Lebesgue's Dominated Convergence Theorem:

\bigskip\noindent\textsc{Proposition 5.2. }\emph{Let }$(X,\Sigma,\mu )$\emph{%
\ be a probability space and set }$\varphi(f)=\int fd\mu$\emph{\ for all }$%
f\in B_{\infty}(\Sigma)$\emph{. Then the set }$\frak{T}=\{\chi_{S}:S\in%
\Sigma\}$\emph{\ is }$\varphi$\emph{-total in }$B_{\infty}(\Sigma)$\emph{.}

\bigskip From this we see that if (11) holds for all measurable
characteristic functions $f$, then $(B_{\infty}(\Sigma),\varphi,\tau)$ is
ergodic by Proposition 4.5(ii), hence $(X,\Sigma,\mu,T)$ is ergodic as
mentioned in Section 2.

Finally, with reference to Proposition 4.6(ii), we note that $\frak{T}^{*}=%
\frak{T}$ for $\frak{T}$ as in Proposition 5.2.

Next we briefly look at von Neumann algebras, as they are well-known
examples of unital $\ast$-algebras. Consider a von Neumann algebra $\frak{M} 
$ and suppose $(\frak{M},\varphi,\tau)$ is a $\ast$-dynamical system. For
example, $\tau$ might be a $\ast$-homomorphism leaving $\varphi$ invariant,
that is to say, $\varphi(\tau(A))=\varphi(A)$ for all $A\in\frak{M}$. Then
the results of Section 4 can be applied directly to $(\frak{M},\varphi,\tau)$%
. As a more explicit (and ergodic) example, we note that $\frak{A}$ in
Example 4.7 is a von Neumann algebra on the Hilbert space $\mathbb{C}^{2}$.
We can also mention that $\tau$ in Example 4.7 is not a homomorphism.

We now describe one suitable choice for the $\varphi$-total sets appearing
in Propositions 4.5 and 4.6. Let $\frak{P}$ be the projections of $\frak{M}$%
. It is known that $\frak{M}$ is the norm closure of the linear span of $%
\frak{P}$, as is mentioned for example on p. 326 of [\textbf{2}]. Since any
state $\varphi$ on $\frak{M}$ is continuous by virtue of being positive, it
follows that $\frak{P}$ is $\varphi$-total in $\frak{M}$. Note also,
regarding Proposition 4.6(ii), that $\frak{P}^{\ast}=\frak{P}$. This is all
very similar to the measure theoretic case in Proposition 5.2, since the
measurable characteristic functions on $X$ are exactly the projections of $%
B_{\infty }(\Sigma)$. This similarity should not be too surprising, since
the theory of von Neumann algebras is often described as ``non-commutative
measure theory'' because of the close analogy with measure theory.

\section*{Appendix}

This Appendix is devoted to the construction of a $\ast$-dynamical system $(%
\frak{A},\varphi,\tau)$ with the property that if $\left\| \tau
(A)-A\right\| _{\varphi}=0$, then $\left\| A-\alpha\right\| _{\varphi}=0$
for some $\alpha\in\mathbb{C}$, but for which the fixed points of the
operator $U$ defined in Proposition 2.4 in terms of some cyclic
representation, form a vector subspace of $\frak{H}$ with dimension greater
than one. This will prove the necessity of a sequence, rather than a single
element, in Definition 2.3, in order for Proposition 2.4 to hold.

First some general considerations. Consider a dense vector subspace $\frak{G}
$ of a Hilbert space $\frak{H}$, and let $\frak{L(H)}$ be the bounded linear
operators $\frak{H\rightarrow H}$. Set 
\begin{equation*}
\frak{A:=}\left\{ A|_{\frak{G}}:A\in\frak{L(H)}\text{, }A\frak{G\subset G}%
\text{ and }A^{\ast}\frak{G\subset G}\right\}
\end{equation*}
where $A|_{\frak{G}}$ denotes the restriction of $A$ to $\frak{G}$. For any $%
A\in\frak{A}$, denote by $\overline{A}$ the (unique) bounded linear
extension of $A$ to $\frak{H}$. Now define 
\begin{equation*}
A^{\ast}:=\overline{A}^{\ast}|_{\frak{G}}
\end{equation*}
for all $A\in\frak{A}$, then it is easily verified that $\frak{A}$ becomes a
unital $\ast$-algebra. (For example, for $A,B\in\frak{A}$ it is clear that $%
AB$ is a bounded linear operator $\frak{G\rightarrow G}$ which therefore has
the extension $\overline{A}.\overline{B}\in\frak{L(H)}$ for which $\overline{%
A}.\overline{B}\frak{G\subset G}$ and $\left( \overline {A}.\overline{B}%
\right) ^{\ast}\frak{G}=\overline{B}^{\ast}\overline{A}^{\ast}\frak{G\subset
G}$ by the definition of $\frak{A}$. Hence $AB\in\frak{A}$, and $\left(
AB\right) ^{\ast}=\left( \overline {A}.\overline{B}\right) ^{\ast}|_{\frak{G}%
}=\left( \overline{B}^{\ast }\overline{A}^{\ast}\right) |_{\frak{G}}=%
\overline{B}^{\ast}\left( \overline{A}^{\ast}|_{\frak{G}}\right) =\overline{B%
}^{\ast}A^{\ast}=B^{\ast }A^{\ast}$. Similarly for the other defining
properties of a unital $\ast $-algebra.) Note that for $A\in\frak{A}$ and $%
x,y\in\frak{G}$ we have 
\begin{equation*}
\left\langle x,Ay\right\rangle =\left\langle x,\overline{A}y\right\rangle
=\left\langle \overline{A}^{\ast}x,y\right\rangle =\left\langle A^{\ast
}x,y\right\rangle \text{.}
\end{equation*}
For a given norm one $\Omega\in\frak{G}$ we define a state $\varphi$ on $%
\frak{A}$ by 
\begin{equation*}
\varphi(A)=\left\langle \Omega,A\Omega\right\rangle \text{.}
\end{equation*}
Next we construct a cyclic representation of $(\frak{A},\varphi)$. Let 
\begin{equation*}
\pi:\frak{A\rightarrow}L(\frak{G}):A\mapsto A
\end{equation*}
then clearly $\pi$ is linear with $\pi(1)=1$ and $\pi(AB)=\pi(A)\pi(B)$.
Note that for any $x,y\in\frak{G}$ we have $(x\otimes y)^{\ast}=y\otimes x$,
hence $(x\otimes y)\frak{G\subset G}$ and $(x\otimes y)^{\ast}\frak{G\subset
G}$, so $(x\otimes y)|_{\frak{G}}\in\frak{A}$. Now, $\pi\left( (x\otimes
\Omega)|_{\frak{G}}\right) \Omega=x\left\langle \Omega,\Omega\right\rangle
=x $, hence $\pi(\frak{A})\Omega=\frak{G}$. Furthermore, $\left\langle
\pi(A)\Omega,\pi(B)\Omega\right\rangle =\left\langle A\Omega,B\Omega
\right\rangle =\left\langle \Omega,A^{\ast}B\Omega\right\rangle =\varphi
(A^{\ast}B)$. Thus $(\frak{G},\pi,\Omega)$ is a cyclic representation of $(%
\frak{A},\varphi)$.

Suppose we have a unitary operator $U:\frak{H}\rightarrow\frak{H}$ such that 
$U\frak{G=G}$ and $U\Omega=\Omega$. Then $U^{\ast}\frak{G}=U^{-1}\frak{G=G}$%
, so $V:=U|_{\frak{G}}\in\frak{A}$, and $V^{\ast}=U^{\ast}|_{\frak{G}}$. It
follows that $VAV^{\ast}\in\frak{A}$ for all $A\in\frak{A}$, hence we can
define a linear function $\tau:\frak{A}\rightarrow\frak{A}$ by 
\begin{equation*}
\tau(A)=VAV^{\ast}\text{.}
\end{equation*}
Clearly $V^{\ast}V=1=VV^{\ast}$, so $\tau(1)=1$ and $\varphi\left(
\tau(A)^{\ast}\tau(A)\right) =\varphi\left( VA^{\ast}AV^{\ast}\right)
=\left\langle U^{\ast}\Omega,A^{\ast}AU^{\ast}\Omega\right\rangle
=\varphi(A^{\ast}A)$, since $U^{\ast}\Omega=U^{-1}\Omega=\Omega$. Therefore $%
(\frak{A},\varphi,\tau)$ is a $\ast$-dynamical system. Note that $U|_{\frak{G%
}}$ satisfies (3), namely $U\pi(A)\Omega=UA\Omega=UAU^{\ast}\Omega=\tau(A)%
\Omega=\pi\left( \tau(A)\right) \Omega$, hence $U$ is the operator which
appears in Proposition 2.4.

Assume $\left\{ x\in \frak{G}:Ux=x\right\} =\mathbb{C}\Omega $. If $\left\|
\tau (A)-A\right\| _{\varphi }=0$, it then follows for $x=\iota (A)$, with $%
\iota $ given by (2), that $\left\| Ux-x\right\| =\left\| \iota \left( \tau
(A)-A\right) \right\| =\left\| \tau (A)-A\right\| _{\varphi }=0$, so $%
x=\alpha \Omega $ for some $\alpha \in \mathbb{C}$. Therefore $\left\|
A-\alpha \right\| _{\varphi }=\left\| \iota (A-\alpha )\right\| =\left\|
x-\alpha \Omega \right\| =0$.

In other words, assuming that the fixed points of $U$ in $\frak{G}$ form the
one-dimensional subspace $\mathbb{C}\Omega$, it follows that $\left\|
\tau(A)-A\right\| _{\varphi}=0$ implies that $\left\| A-\alpha\right\|
_{\varphi}=0$ for some $\alpha\in\mathbb{C}$.

It remains to construct an example of a $U$ with all the properties
mentioned above, whose fixed point space in $\frak{H}$ has dimension greater
than one. The following example was constructed by L. Zsid\'{o}:

Let $\frak{H}$ be a separable Hilbert space with an orthonormal basis of the
form 
\begin{equation*}
\{\Omega \,,\,y\}\cup \left\{ u_{k}\,:\,k\in \mathbb{Z}\right\} 
\end{equation*}
(that is to say, this is a total orthonormal set in $\frak{H}$) and define
the linear operator $U:\frak{H}\longrightarrow \frak{H}$ by 
\begin{align*}
U\Omega & =\Omega \,, \\
Uy& =y\,, \\
Uu_{k}& =u_{k+1}\,,\quad k\in \mathbb{Z}\text{.}
\end{align*}
Since $U$ is a surjective isometry, it is unitary. Let $\frak{G}$ be the
linear span of 
\begin{equation*}
\{\Omega \}\cup \left\{ y+u_{k}:\,k\in \mathbb{Z}\right\} \text{.}
\end{equation*}
Then $U\frak{G}=\frak{G}$. Furthermore, $\frak{G}$ is dense in $\frak{H}\,$.
Indeed, 
\begin{equation*}
\Vert y-\frac{1}{n}\sum_{k=1}^{n}(y+u_{k})\Vert =\frac{1}{n}\Vert
\sum_{k=1}^{n}u_{k}\Vert =\frac{1}{\sqrt{n}}\longrightarrow 0
\end{equation*}
implies that $y\in \overline{\frak{G}}\,$, the closure of $\frak{G}$, hence
also 
\begin{equation*}
u_{k}=(y+u_{k})-y\in \overline{\frak{G}}\,
\end{equation*}
for $k\in \mathbb{Z}$.

Next we show that 
\begin{equation*}
\{x\in\frak{G}\,:\,Ux=x\}=\mathbb{C}\Omega.
\end{equation*}
If $\alpha\Omega+\sum\limits_{k=-n}^{n}\beta_{k}(y+u_{k})\in\frak{G}$ is
left fixed by $U$, then 
\begin{equation*}
\alpha\Omega+\sum\limits_{k=-n}^{n}\beta_{k}y+\sum\limits_{k=-n}^{n}\beta
_{k}u_{k+1}=\alpha\Omega+\sum\limits_{k=-n}^{n}\beta_{k}y+\sum%
\limits_{k=-n}^{n}\beta_{k}u_{k}
\end{equation*}
and it follows that $\beta_{-n}=0$, and that $\beta_{k+1}=\beta_{k}$ for $%
k=-n,...,n-1$. Thus 
\begin{equation*}
\alpha\Omega+\sum\limits_{k=-n}^{n}\beta_{k}(y+u_{k})=\alpha\Omega\,.
\end{equation*}

On the other hand, 
\begin{equation*}
\{x\in\frak{H}\,:\,Ux=x\}
\end{equation*}
clearly contains the two-dimensional vector space spanned by $\Omega$ and $%
y\,$.

\bigskip\noindent\emph{Acknowledgements. }We would like to thank the
National Research Foundation and the Mellon Foundation Mentoring Programme
for financial support. We are also grateful to L\'{a}szl\'{o} Zsid\'{o} for
providing us with the example used in the Appendix.

\end{document}